# Time-Energy Optimal Control of a Mobile Robot Using Direct Numerical Method

Mohamad Shahab, Amar Khoukhi, and Fouad Al-Sunni

*Abstract*—Optimal control of a mobile robot system is formulated. Multiobjective criteria of time and energy is employed. The optimal control problem is formulated as a nonlinear programming problem (NLP). The problem is solved using the direct method of numerical optimal control. This setting showed great flexibility in incorporating different information relating to the problem, namely physical constraints and nonlinear dynamics of the system. System inputs are considered as optimization variables, along with sampling periods of the applied inputs being optimization variables as well. Different scenarios on objectives of the problem are implemented and investigated. Interesting results are found in terms of complying with the expected behavior of a mobile robot system.

## I. Introduction

RESEARCH on mobile robot systems is always representing an interesting field of investigation. In the same time, challenges in dealing with these systems also are increasing. This work investigates the problem of finding a solution for optimal control of a nonholonomic mobile robot. A multiple objectives criterion (time & energy) is considered here. Time and energy have always been adversaries in terms incorporating both quantities in one setting. optimal control of robotic systems are somehow difficult, mainly because of nonlinearities involved. So, different approaches were tackled in the literature.

### A. Numerical Optimal Control: Direct Method

Generally, numerical solutions for optimal control problems are needed when handling nonlinear systems. Unlike linear systems, analytic solutions are hard or even non-existent. Numerical solutions of optimal control problems are divided into two main categories, namely, direct methods and indirect methods. Indirect methods are the ones that adhere to the 1st-order optimality conditions via Pontryagin's Minimum Principle or Variational principle [1]. The Hamiltonian systems would be solved numerically as a Two Boundary-value problem.

The great downside of indirect methods is that they cannot handle discontinuities in the problem. Discontinuities can come from the introduction of constraints to the optimization problem or the high nonlinearities in the system. So, direct methods come to replace indirect ones to solve for more complex optimal control problems. Here in this work, direct methods are considered. The basic principle of direct methods is to discretize the problem and solve it using well-established Nonlinear Programming (NLP) techniques.

Direct methods can be applied through different variants. Here, Sequential (Single Shooting) approach will be used [3]. Other approaches are more investigated in [3]. The basic idea behind direct sequential method of numerically finding the solution of an optimal control problem is to have the system input variables discretized at time instants. Then, system states can be computed throughout the time horizon via any ordinary differential equation solvers considering the time-discretized input.

### B. General Optimal Control Problem

So, a general optimal control problem can be formulated as,

$$\min_u H(x(t_f)) + \int_{t_0}^{t_f} \mathcal{L}(x(t), u(t), t) \, dt$$

$$s.t.$$
$$\dot{x} = f(x, u), \ x(t_0) = x_0 \qquad (1)$$
$$\mathbf{g}(x, u) \leq 0$$
$$\mathbf{h}(x, u) = 0$$

The above optimal control generic problem has a final state objective of $H(x(t_f))$ and the Lagrangian objective of $\mathcal{L}(x(t), u(t), t)$. System behavior is governed by the nonlinear dynamic system in $\dot{x} = f(x, u)$ with an initial condition of $x(t_0) = x_0$. The performance is restricted by collection of inequality and equality constraints of $\mathbf{g}(x, u) \leq 0$ & $\mathbf{h}(x, u) = 0$, respectively.

This paper is organized as follows. Section II will give the system model of the nonholonomic mobile robot in relation with the algorithm at hand. In section III, the Multi-objective optimal control formulation of the problem will be explained in detail. Implementation will be provided in IV with discussion about the results. Section V will conclude the paper.

Mohamad Shahab is a research assistant with the Department of Systems Engineering, King Fahd University of Petroleum and Minerals, Dhahran, Saudi Arabia. (e-mail: moh_shahab@hotmail.com).

Amar Khoukhi is assistant professor with the Department of Systems Engineering, King Fahd University of Petroleum and Minerals, Dhahran, Saudi Arabia. (e-mail: amar@kfupm.edu.sa).

Fouad Alsunni is a professor with the Department of Systems Engineering, King Fahd University of Petroleum and Minerals, Dhahran, Saudi Arabia. (e-mail: alsunni@kfupm.edu.sa).

## II. MOBILE ROBOT SYSTEM MODEL

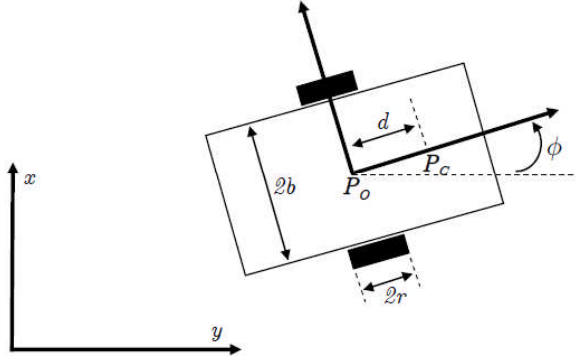

Fig. 1. Nonholonomic Mobile Robot

This work investigates the solution of a time-energy optimal control problem for a nonholonomic mobile robot. For a wheeled mobile robot as in figure 1, configuration can be defined by both position $(x,y)$ and orientation $\phi$. Systems dynamics can be described by the nonlinear system of $\dot{\mathbf{x}} = f(\mathbf{x},\tau)$, with $\mathbf{x} = \begin{bmatrix} x & y & \phi & \theta_R & \theta_L & v_R & v_L \end{bmatrix}^T$, with $(\theta_R, \theta_L, v_R, v_L)$ being right and left wheels angular positions and velocities, respectively. We have $f(\mathbf{x},\tau)$ as

$$f(\mathbf{x},\tau) = \begin{bmatrix} 0 & S \\ 0 & -M^{-1}V \end{bmatrix}\mathbf{x} + \begin{bmatrix} 0 \\ M^{-1} \end{bmatrix}\tau \quad (2)$$

Matrices S, M, and V are corresponding to kinematic relation matrix, inertia matrix, and coriolis matrix, respectively (see [7] for more details). The above dynamics can be linearized [6] by feedback of $\tau = Mu + V\nu$, with $\nu = \begin{bmatrix} v_R & v_L \end{bmatrix}^T$ and the auxiliary input of $u$. So, the systems can be formed as,

$$\dot{\mathbf{x}} = f(\mathbf{x}) + g(\mathbf{x})u = \begin{bmatrix} S\nu \\ 0 \end{bmatrix} + \begin{bmatrix} 0 \\ I \end{bmatrix}u \quad (3)$$

The above system can be viewed as two parts, namely, the kinematic model of $S\nu$, and the auxiliary acceleration inputs of $\dot{\nu} = u$.

In order to solve the problem in a Nonlinear Programming (NLP) framework, input variable should be discretized through time instants. Let us consider input

$$u(t) = u(kT_s) = u(k) = \text{constant}, \quad \text{for } kT_s \leq t < kT_s + T_s \quad (4)$$

with $t$ being the time independent variable, $k$ as being the time index in discrete domain, and $T_s$ being the sampling period. In principle, system inputs can be discretized in any other different fashion. Here, a piecewise constant (zero-order) input is considered. A piecewise linear (first-order) can also be employed.

During the optimization algorithm iterations, system states should be computed. State values are used for the objective function and constraints throughout the whole optimization time horizon. Discrete-time model can be computed using the Taylor expansion as [4],

$$\mathbf{x}(k+1) = \mathbf{x}(k) + \sum_{\ell=1}^{\ell_{\max}} D^{[\ell]}(\mathbf{x}(k), u(k)) \frac{T_s^\ell}{\ell!} \quad (5)$$

With $D^{[\ell]}(\cdot,\cdot)$ is found recursively by,

$$D^{[1]}(\mathbf{x},u) = f(\mathbf{x}) + g(\mathbf{x})u$$
$$D^{[\ell+1]}(\mathbf{x},u) = \frac{\partial D^{[\ell]}(\mathbf{x},u)}{\partial \mathbf{x}} \cdot (f(\mathbf{x}) + g(\mathbf{x})u) \quad (6)$$

Observe that *exact* discretization can be computed when $\ell_{\max} = \infty$. Throughout this work, $\ell_{\max}$ is chosen to be 2 which is resembling a 2nd-order Taylor expansion.

The discrete-time model of the mobile robot having both accelerations and sampling-period as system inputs can be put as:

$$\mathbf{x}(k+1) = f_D(\mathbf{x}(k), u(k), T_s(k)) \quad (7)$$
$$= \mathbf{A}(\mathbf{x}(k), T_s(k))\mathbf{x}(k) + \mathbf{B}(\mathbf{x}(k), T_s(k))u(k)$$

with,

$$\mathbf{A}(\mathbf{x}(k), T_s(k)) = \begin{bmatrix} \mathbf{I}_{5\times5} & \overline{\mathbf{A}}(\mathbf{x}(k), T_s(k)) \\ \mathbf{0}_{2\times5} & \mathbf{I}_{2\times2} \end{bmatrix}$$

And,

$$\overline{\mathbf{A}}(\mathbf{x}(k), T_s(k)) = \begin{bmatrix} (\frac{T_s r}{2}\cos\phi(k) - \frac{T_s^2 r^2}{8b}(v_R(k)+v_L(k))\sin\phi(k)) & (\frac{T_s r}{2}\cos\phi(k) + \frac{T_s^2 r^2}{8b}(v_R(k)+v_L(k))\sin\phi(k)) \\ (\frac{T_s r}{2}\sin\phi(k) + \frac{T_s^2 r^2}{8b}(v_R(k)+v_L(k))\cos\phi(k)) & (\frac{T_s r}{2}\sin\phi(k) - \frac{T_s^2 r^2}{8b}(v_R(k)+v_L(k))\cos\phi(k)) \\ \frac{T_s(k)r}{2b} & -\frac{T_s(k)r}{2b} \\ T_s(k) & 0 \\ 0 & T_s(k) \end{bmatrix}$$

Also,

$$\mathbf{B}(\mathbf{x}(k), T_s(k)) = \begin{bmatrix} \frac{T_s^2(k)r}{4}\cos\phi(k) & \frac{T_s^2(k)r}{4}\cos\phi(k) \\ \frac{T_s^2(k)r}{4}\sin\phi(k) & \frac{T_s^2(k)r}{4}\sin\phi(k) \\ \frac{T_s^2(k)r}{4b} & -\frac{T_s^2(k)r}{4b} \\ \frac{T_s^2(k)}{2} & 0 \\ 0 & \frac{T_s^2(k)}{2} \\ T_s & 0 \\ 0 & T_s \end{bmatrix}$$

With $r$ being the radius of the wheels, and $b$ being the radius of the robot body. The above discrete model can be imagined with 'two' input sets of variables. The first set resembles the discretized acceleration inputs of $u(k) = \begin{bmatrix} u_R(k) & u_L(k) \end{bmatrix}^T$. The other set consist of the sampling period as a control variable [5] in the optimization problem $T_s(k)$. So, optimization problem control variables are $u(k), T_s(k)$.

## III. TIME-ENERGY FORMULATION

### A. Optimal Control Problem Formulation

Having the problem to be set for a Nonlinear Programming (NLP) setting, considering the objective function of

$$\min H(\mathbf{x}(t_f)) + \int_{t_0}^{t_f} \mathcal{L}(\mathbf{x}(t), u(t), t) \, dt \qquad (8)$$

The objective can be reformulated into the Mayer form. The Mayer form gives big advantage in terms of the computer algorithm. The new objective can be put as,

$$\min H(\mathbf{x}(t_f)) + z(t_f) \qquad (9)$$

With the Lagrangian information is embedded into a dummy state variable $z(t)$,

$$\dot{z} = \mathcal{L}(\mathbf{x}(t), u(t), t) \qquad (10)$$

So, with review of original problem in (1), the problem can be now formulated as,

$$\min_{\{u, T_s\}} H(\mathbf{x}(N)) + z(N) \qquad (11)$$

s.t.

$z(k+1) = z(k) + T_s(k) \cdot \{\mathcal{L}(\mathbf{x}(k), u(k), T_s(k))\}, z(0) = 0$
$\mathbf{x}(k+1) = f_D(\mathbf{x}(k), u(k), T_s(k)), \ k = 1, \ldots, N, \mathbf{x}(0) = \mathbf{x}_0$
$\mathbf{g}(\mathbf{x}, u, T_s) \leq 0$
$\mathbf{h}(\mathbf{x}, u, T_s) = 0$

Number of time instants across time horizon is put as $N$. Optimal Control performance measures would be embedded inside the Lagrangian, $\mathcal{L}(\mathbf{x}, u, T_s)$.

So, the above optimization problem is now formulated in a Nonlinear Programming (NLP) setting. Optimization control variables are $[u_R(k), u_L(k)], T_s(k), k = 1 \ldots N$. So, the number of optimization variable are $3 \times N$. If time-sampling period is not considered as a variable, $2 \times N$ would be the number of optimization variables. This NLP problem can be solved by any state-of-the-art algorithm available.

### B. Time-energy Optimal Control

Here, energy and time are considered for optimization. For wheeled mobile robot, energy can be viewed in two sets: input energy, and kinetic energy. Input energy can be computed from the wheels torques spent throughout the process. For actual physical system, some constraints should be restricting the process. Let us consider our optimal control problem of,

$$\min_{\{u, T_s\}} H(\mathbf{x}(N)) + z(N)$$

s.t.

$z(k+1) = z(k) + T_s(k) \cdot \{\mathcal{L}(\mathbf{x}(k), u(k), T_s(k))\}, z(0) = 0$
$\mathbf{x}(k+1) = f_D(\mathbf{x}(k), u(k), T_s(k)), \ k = 1, \ldots, N, \mathbf{x}(0) = \mathbf{x}_0$
$\mathbf{g}(\mathbf{x}, u, T_s) \leq 0$
$\mathbf{h}(\mathbf{x}, u, T_s) = 0$

#### 1) Cost functions

Let us have the Lagrangian for the problem be,

$$\mathcal{L}(\mathbf{x}, u, T_s) = E_{IN} + E_{KE} + \beta \qquad (12)$$

We have $E_{IN}$ as the cost for energy spent by the torques of the wheels, and $E_{KE}$ as the cost for kinetic energy spent by robot body, and $\beta$ as the weight on time. Quadratic measures are considered for the energies.

$$E_{IN} = \tau^T R \tau \qquad (13)$$

With $R$ as a diagonal matrix to weight the input energy. The torque $\tau$ is calculated as the feedback linearization function (function of the acceleration inputs and angular velocities of right and left wheels),

$$\tau = Mu + V\nu \qquad (14)$$

For the kinetic energy cost, we can put it as,

$$E_{KE} = \nu^T P \nu \qquad (15)$$

With $P$ as a diagonal matrix to weight the kinetic energy. Velocities considered are $\nu = \begin{bmatrix} v_R & v_L \end{bmatrix}^T$. However, another cost measure can be calculated as,

$$E_{KE} = \mathbf{v}^T P \mathbf{v} \qquad (16)$$

With,

$$\mathbf{v} = \Omega^{-1} \nu \qquad (17)$$

With $\mathbf{v} = \begin{bmatrix} v & \omega \end{bmatrix}^T$ as the robot body velocities of linear velocity, $v$, and angular velocity, $\omega$. The above transformation is between the wheel velocities and the robot body velocities which are related via the geometric information in,

$$\Omega = \begin{bmatrix} \dfrac{1}{r} & \dfrac{b}{r} \\ \dfrac{1}{r} & -\dfrac{b}{r} \end{bmatrix} \qquad (18)$$

For the final state cost of $H(\mathbf{x}(N))$, let have it be zero for this work. So, the objective of the problem be,

$$\min_{\{u, T_s\}} z(N) \qquad (19)$$

#### 2) Constraints

For the problem in (11), constraints are considered to restrict the process according to physical information available. With the sampling period being an input, upper and lower bounds are enforced. More relevant bounds are the bounds on the torques provided by the electrical motors to the wheels. To assure the robot system to reach a desired configuration, final states are constrained also.

Bounds on sampling period can be put as $T_{\min} \leq T_s(k) \leq T_{\max}$. So, two constraints put in standard form are,

$$\begin{aligned} [g_1] &= T_s(k) - T_{\max} \leq 0, k = 1 \ldots N \\ [g_2] &= T_{\min} - T_s(k) \leq 0, k = 1 \ldots N \end{aligned} \qquad (20)$$

With torques computed via (14), and bounds of $\tau_{\min} \leq \tau \leq \tau_{\max}$, constraints can be formulated as,

$$[g_3] = Mu(k) + V\nu(k) - \tau_{\max} \leq 0, k = 1...N$$
$$[g_4] = \tau_{\min} - Mu(k) + V\nu(k) \leq 0, k = 1...N \quad (21)$$

To ensure the robot to reach its destination in both position and orientation, constraints on the final state for the system should be employed. A desired final configuration is,

$$x(N), y(N), \phi(N) = x_F, y_F, \phi_F$$

A more tolerant yet acceptable formulation, constraints can be put as,

$$\begin{aligned} g_5 &= |x(N) - x_F| - \varepsilon \leq 0 \\ g_6 &= |y(N) - y_F| - \varepsilon \leq 0 \\ g_7 &= |\phi(N) - \phi_F| - \varepsilon_\phi \leq 0 \end{aligned} \quad (22)$$

With $\varepsilon, \varepsilon_\phi$ are tolerances in Cartesian dimensions and in angle dimension, respectively.

Also, in order to have the mobile robot to reach its destination regulated to rest, final velocities constraint is,

$$g_8 = |\mathbf{v}(N)| - \varepsilon_v \leq 0 \quad (23)$$

With final body velocities $\mathbf{v}$ are computed from (17) restricted with tolerance of $\varepsilon_v$. So, the complete set of inequality constraints are,

$$\mathbf{g}(\mathbf{x}, u, T_s) \leq 0 \quad (24)$$

With

$$\mathbf{g}(\mathbf{x}, u, T_s) = \begin{bmatrix} [g_1]_{\forall k} & [g_2]_{\forall k} & [g_3]_{\forall k} & [g_4]_{\forall k} & g_5 & g_6 & g_7 & g_8 \end{bmatrix}$$

So, we can have up to $4 \times N + 4$ inequality constraints for the optimization problem and no equality constraints. Observe that constraints 5 through 8 are equality constraints when not tolerating, i.e. tolerances values are all zeros. Further, other constraints can be imagined for the process. Bounds on accelerations, like $u_{\min} \leq u(k) \leq u_{\max}$ and velocities, like $\nu_{\min} \leq \nu(k) \leq \nu_{\max}$, can all be enforced.

The complete Nonlinear Programming (NLP) problem of time-energy optimal control of a mobile robot can be viewed in (25), with corresponding objectives (12) and constraints in (24).

$$\begin{aligned} &\min_{\{u, T_s\}} z(N) \\ &s.t. \\ &z(k+1) = z(k) + T_s(k) \cdot \{\mathcal{L}(\mathbf{x}(k), u(k), T_s(k))\}, z(0) = 0 \quad (25) \\ &\mathbf{x}(k+1) = f_D(\mathbf{x}(k), u(k), T_s(k)), \ k = 1,...,N, \mathbf{x}(0) = \mathbf{x}_0 \\ &\mathbf{g}(\mathbf{x}, u, T_s) \leq 0 \end{aligned}$$

## IV. IMPLEMENTATION AND DISCUSSION

Here, solution for the time-energy optimal control problem would be explored. Energy and time costs will go head to head for having a solution to the problem. Under MATLAB, the NLP problem in (25) is solved utilizing Sequential Quadratic Programming (SQP) algorithm with Hessian computation using quasi-Newton (BFGS) approximation [2]. Algorithm worked with termination tolerance of $10^{-9}$.

Consider the costs in (12), let's assume

$$R = P = \alpha \mathbf{I}_{2\times 2} \quad (26)$$

to have an identical weight $\alpha$ on both kinetic and input energies. Also, weight on process time $\beta$.

To quantify the performance of the solution, let have the optimal solution of $u^*(k), T_s^*(k)$, $k = 1...N$. A measure for the optimal total energy can be computed as,

$$E_{optimal} = \sum_{k=1}^{N} T_s^*(k) \cdot \left[ k_\tau \cdot \tau^*(k)^T \tau^*(k) + k_m \cdot \mathbf{v}^*(k)^T \mathbf{v}^*(k) \right] \quad (27)$$

Energy in (27) is calculated with appropriate gains, $k_m, k_\tau$, for kinetic energy, and electric input energy, respectively. Also, the optimal final time can be computed as,

$$t_f = \sum_{k=1}^{N} T_s^*(k) \quad (28)$$

Now, the NLP problem to be solved is formulated in (25). Let us have $N = 40$. The mobile robot system has the detailed model information in [7] with all the required geometric and physical parameters values of matrices S, M, and V.

For the constraints information in (24), let have the values of $T_{\min} = 0.01 \text{sec}$, $T_{\max} = 2 \text{sec}$ & $\tau_{\min} = -1$, $\tau_{\max} = +1$. Constraints tolerances are put as $\varepsilon = 10^{-3} \text{m}, \varepsilon_\phi = 10^{-7} \text{rad}, \varepsilon_v = 10^{-9} \text{m/sec}$.

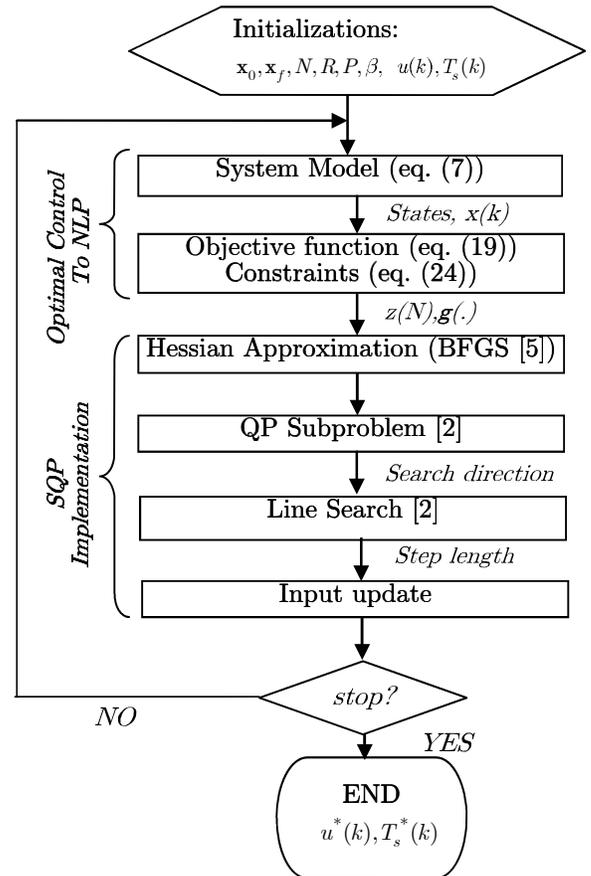

Fig. 2. Overview Flowchart of the algorithm

The multiobjective problem is solved for multiple weights on time and energy. A weight of $\beta = 0$ means an only energy minimization problem. A weight of $\alpha = 0$ means an only time minimization problem.

### A. Example 1

A first example of solutions for the work is to have initial and final configurations of:

$$\begin{bmatrix} x_0 \\ y_0 \\ \phi_0 \end{bmatrix} = \begin{bmatrix} 7 \\ 1 \\ \frac{\pi}{2} \end{bmatrix}, \quad \begin{bmatrix} x_F \\ y_F \\ \phi_F \end{bmatrix} = \begin{bmatrix} -7 \\ -1 \\ -\frac{\pi}{2} \end{bmatrix}.$$

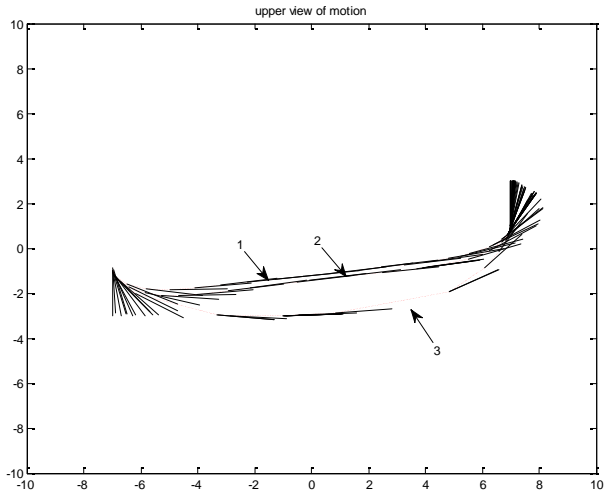

Fig. 3. optimal trajectories for example 1, 1) time & energy, 2) pure energy, and 3) pure time. Solid black lines represent the orientations.

For these configurations, three solutions would be explored, namely, 1) time & energy minimization, 2) pure energy minimization, and 3) pure time minimization. Table 1 shows the results of these weightings of the problem. Total energy and total time is computed using (27), and (28). Figure 3 shows the upper view for the 3 optimal motions. To see it more clearly, the time-minimization solution gave a bang-bang-like solution, as in figure 4 showing the optimal torques. Figure 5 shows the optimal time periods for torques to be applied in.

### B. Example 2

In this example, 3 different initial and final configurations is put. Energy and time minimization is employed with $\alpha = 10, \beta = 10$. The three initial and final configurations are:

1) $\mathbf{x}_0 = [7, 7, 0]$, $\mathbf{x}_F = [-7, -7, \pi]$,
2) $\mathbf{x}_0 = [-7, 0, \pi]$, $\mathbf{x}_F = [7, 0, 0]$, and
3) $\mathbf{x}_0 = [-7, 0, \frac{\pi}{4}]$, $\mathbf{x}_F = [7, 0, \frac{3\pi}{4}]$.

Figure 6 shows the 3 upper views of the motions. Figure 7 shows the errors in x- & y-positions and orientation for each of the three cases throughout the optimized motions.

### C. Example 3

For weights of $\alpha$ as in (26) and $\beta$, values tried are $\alpha \in \{0, 5, 10, 15, 20\}$, $\beta \in \{0, 5, 10, 15, 20\}$. So, different combinations of weights would be enforced into time and energy, namely 25 combinations. These combinations would provide 25 different optimal solutions and trajectories for the mobile robot. To implement the multiple weightings, a generic initial and desired final positions and orientations for the robot is selected in,

$$\begin{bmatrix} x_0 \\ y_0 \\ \phi_0 \end{bmatrix} = \begin{bmatrix} 0 \\ 8 \\ -\frac{\pi}{2} \end{bmatrix}, \quad \begin{bmatrix} x_F \\ y_F \\ \phi_F \end{bmatrix} = \begin{bmatrix} 0 \\ -8 \\ \frac{\pi}{2} \end{bmatrix}$$

We will have 25 different optimal trajectories. Each trajectory corresponds to different combination of weights on time and energy.

Now, from all optimal solutions, figure 8 shows optimal final times versus weights applied on time. Each solid line correspond to each of weights on energy. Figure 9 shows optimal energies versus weights applied on energy. Each solid line correspond to each of weights on time. To have the overall picture of the results of this multiobjective optimization, figure 10 shows normalized surfaces of optimal final time and optimal total energy versus the weights on time and energy.

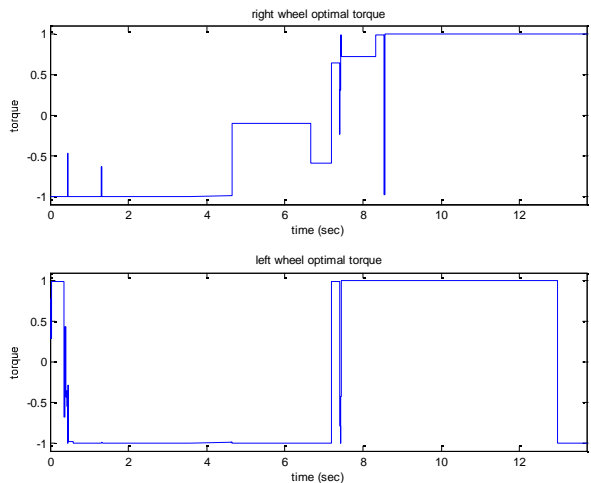

Fig. 4. Right and left wheel time-optimal torques

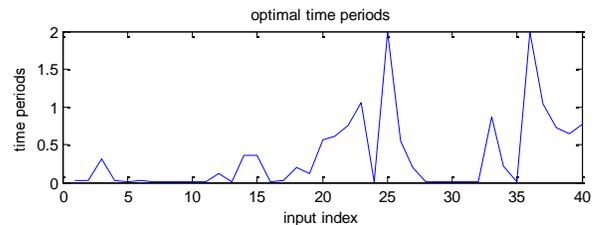

Fig. 5. Optimal time periods for torques to be applied

TABLE I
OPTIMAL VALUES FOR EXAMPLE 1

| Weights | Total Time (sec) | Total Energy (joule) |
|---|---|---|
| $\alpha = 20, \beta = 20$ | 23 | 16.1 |
| $\alpha = 20, \beta = 0$ | 40 | 8.53 |
| $\alpha = 0, \beta = 20$ | 13.73 | 52.85 |

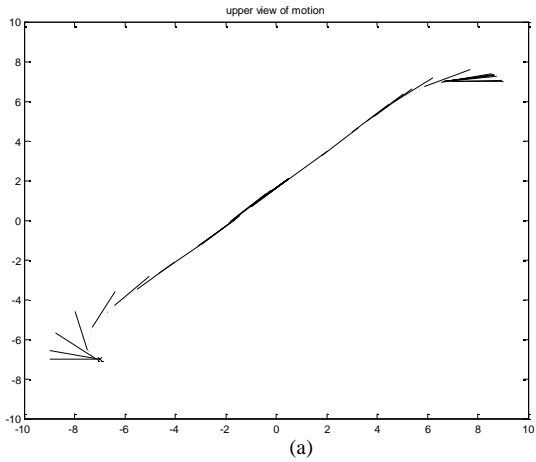

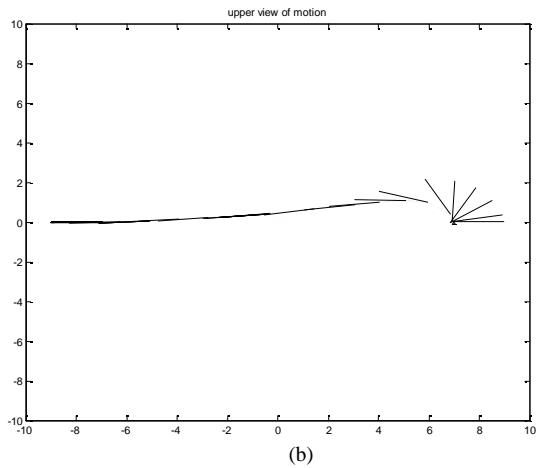

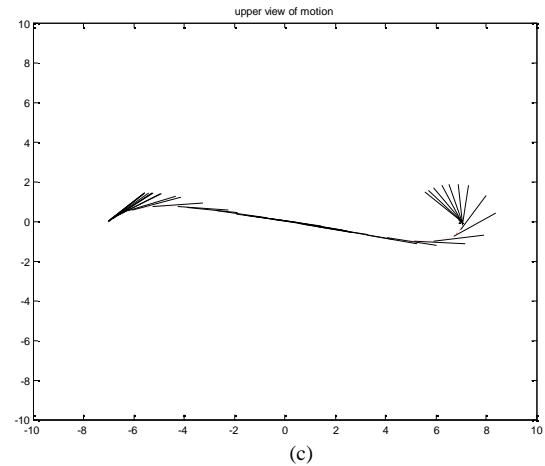

Fig. 6. Upper view of optimal motions of 3 different initial and final configurations. Solid black lines represent the orientations.

### D. Comments and Discussion

From figures, you can see clearly the behavior of the problem. In most of the cases, logical results appear. Putting more weight on energy gives higher minimum time values. In the other hand, more weight on time gives higher minimum energy values. We have to put in mind that this optimization is applied into a mobile robot system with acceleration inputs, i.e. the kinodynamic problem rather than the kinematic problem. You can see that geometric properties of each trajectory are different according to the multiple scenarios of weighing on time and energy.

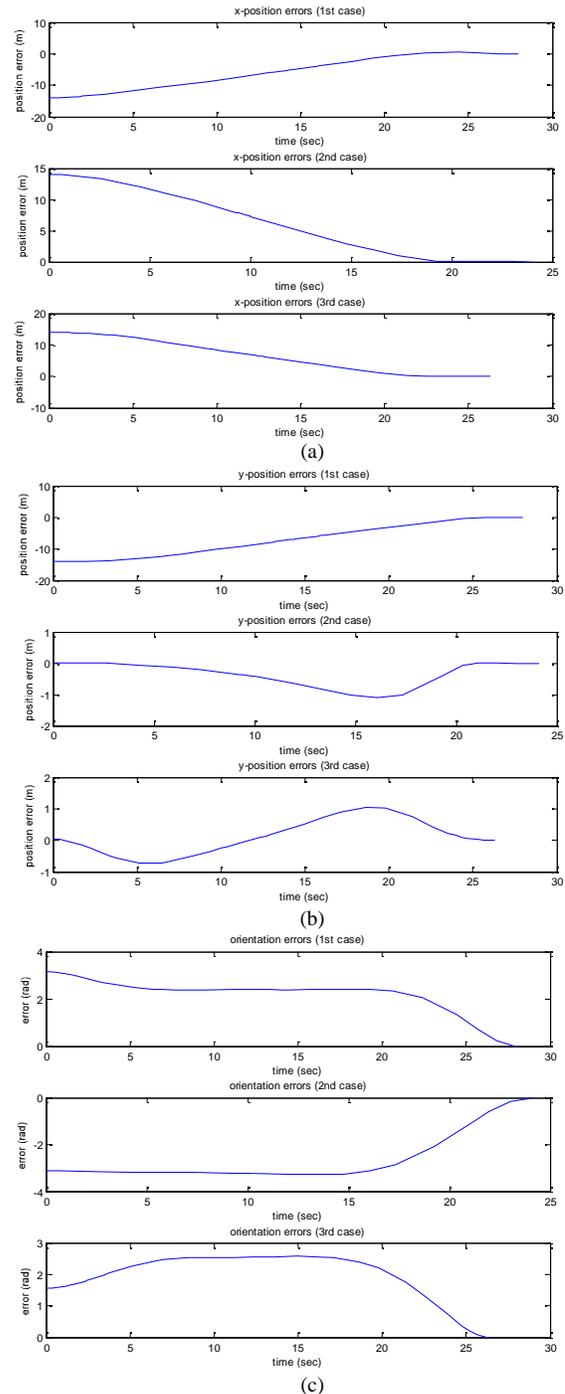

Fig. 7. Errors in the 3 different cases in example 2. Errors in x-positions (a), errors in y-positions (b), and errors in orientations (c).

First example demonstrates the two extreme problems of pure time minimization and pure energy minimization. According to table 1, the difference between the optimal time in case 3 and total time in case 2 shows the sensitivity of the problem for time optimization. In the same time, the incorporation of input sampling period as an optimization variable gave us great flexibility.

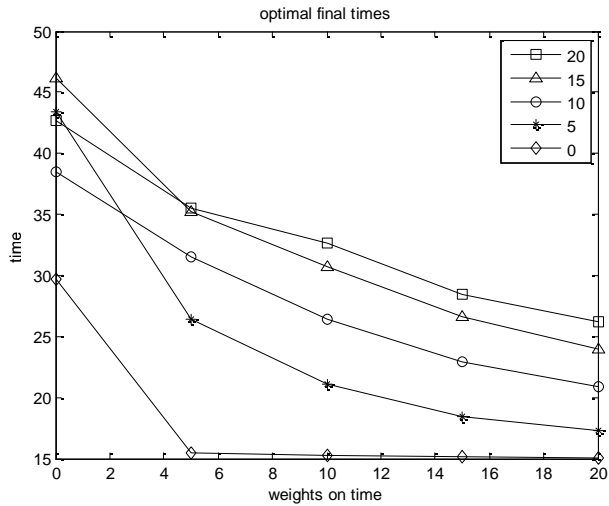

Fig. 8. Optimal final times. Each solid line corresponds to different weight on energy. The values of the weights on energy is shown in the legend.

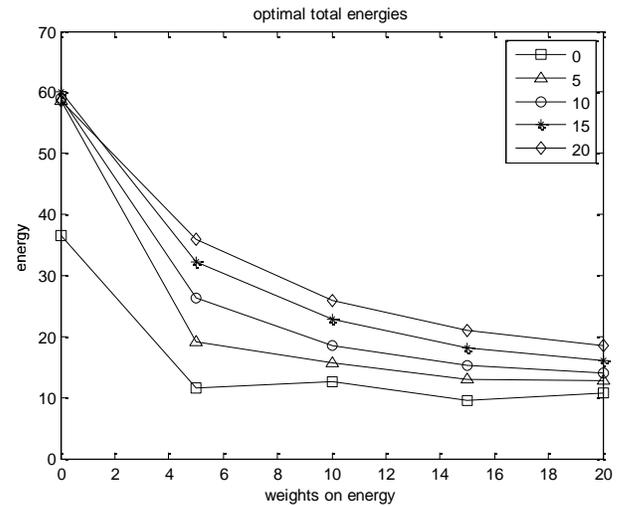

Fig. 9. Optimal total energy. Each solid line corresponds to different weight on time. The values of the weights on time is shown in the legend.

Second example shows us the applicability of the suggested solution for different maneuver scenarios. Putting in mind the regulation to zero velocity at final configurations, attained solutions give a great advantage in many real life applications.

The surfaces in the figure 8 shows a clear picture about the behavior. With defined final time and total energy in (27) and (28) respectively, the result gives insight about how best to optimize the process. Compromise between energy and time can be investigated. From observing the intersection of the two surfaces in figure 10, weighing energy with value around 5 gives us interesting results. As weighing energy with more than five will make optimal time values to rise.

## V. CONCLUSION

Future investigation can be taken to involve more complexities to the problem. With the powerful tools of NLP algorithms, more difficult constraints and objectives can be tested easily. A next step is to put obstacles in the robot environment to make the optimal control problem to handle them. Caution, however, should be taken into consideration as nonconvexities could appear.

Throughout this paper, a multiobjective optimal control problem is studied. Optimal control of a mobile robot system is formulated. The optimal control problem is designed as a Nonlinear Programming problem. This corresponds to the direct approach for numerically solving the optimal control problem. This setting showed great flexibility in incorporating different information relating to the problem, namely physical constraints and nonlinear dynamics of the system. Standard optimal control method lacked this flexibility. Different scenarios on objectives of the problem is investigated. Multiple objectives are put into the problem. Both time and energy are employed as a criteria. Interesting results are found in terms of complying with the expected behavior of a mobile robot system.

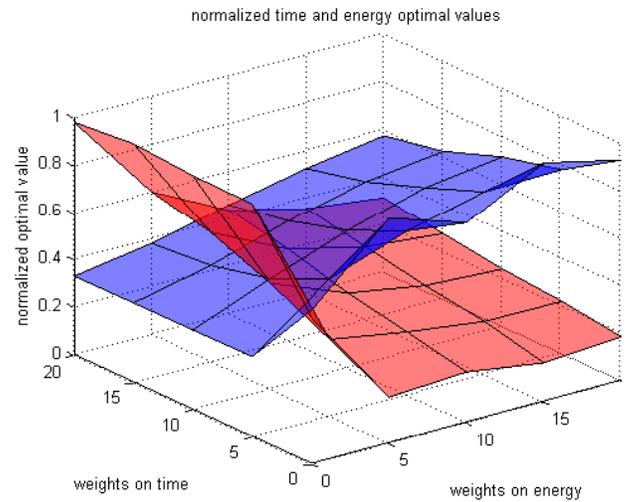

Fig. 10. Normalized Optimal surfaces of time (blue) and energy (red).

ACKNOWLEDGMENT

The authors would like to thank King Fahd University of Petroleum and Minerals for its support for this work.